\documentclass[12pt]{article}

\usepackage{amsfonts} 

\begin{document}

\renewcommand{\baselinestretch}{1.2} 
\normalsize

\newtheorem{lemma}{Lemma}[section]
\newtheorem{thm}[lemma]{Theorem}
\newtheorem{cor}[lemma]{Corollary}
\newtheorem{voorb}[lemma]{Example}
\newtheorem{rem}[lemma]{Remark}
\newtheorem{prop}[lemma]{Proposition}

\newcommand{\Ni}{{\mathbb{N}}}
\newcommand{\Ri}{{\mathbb{R}}}
\newcommand{\Ci}{{\mathbb{C}}}
\newcommand{\Ti}{{\mathbb{T}}}
\newcommand{\Zi}{{\mathbb{Z}}}

\newcommand{\PPi}{{\mathbb{P}}}
\newcommand{\Di}{{\mathbb{D}}}

\newcommand{\proof}{\mbox{\bf Proof.} \hspace{3pt}}
\newcommand{\remark}{\mbox{\bf Remark} \hspace{3pt}}

\newcommand{\ruimte}{\vskip10.0pt plus 4.0pt minus 6.0pt}

\newcommand{\ad}{{\mathop{\rm ad}}}
\newcommand{\Ad}{\mathop{\rm Ad}}
\newcommand{\arccot}{\mathop{\rm arccot}}
\newcommand{\codim}{\mathop{\rm codim}}
\newcommand{\RRe}{\mathop{\rm Re}}
\newcommand{\IIm}{\mathop{\rm Im}}
\newcommand{\HS}{{\mathop{\rm HS}}}
\newcommand{\Tr}{{\mathop{\rm Tr}}}
\newcommand{\supp}{\mathop{\rm supp}}
\newcommand{\sgn}{\mathop{\rm sgn}}
\newcommand{\esssup}{\mathop{\rm ess\,sup}}
\newcommand{\essinf}{\mathop{\rm ess\,inf}}
\newcommand{\liminferior}{\mathop{\rm lim\, inf}}
\newcommand{\Leibniz}{\mathop{\rm Leibniz}}
\newcommand{\lcm}{\mathop{\rm lcm}}
\newcommand{\mod}{\mathop{\rm mod}}
\newcommand{\spann}{\mathop{\rm span}}
\newcommand{\ubar}{\underline{\;}}

\newfont{\fontcmrten}{cmr10}
\newcommand{\slbrl}{\mbox{\fontcmrten (}}
\newcommand{\slbrr}{\mbox{\fontcmrten )}}

\newfont{\fontcmsyseventeen}{cmsy10 scaled\magstep3}
\newcommand{\biggertimes}{{\mbox{{\fontcmsyseventeen \symbol{'002}}}}}

\newenvironment{remarkn}{\begin{rem} \rm}{\end{rem}}
\newenvironment{exam}{\begin{voorb} \rm}{\end{voorb}}

{\bf \Large Subordinated discrete semigroups} 

{\bf \Large of operators}    

\bigskip 

{\bf Nick Dungey} 

\bigskip

\noindent {\bf Abstract.}   
Given a power-bounded linear operator $T$ in a Banach space 
and a probability $F$ on the non-negative integers, 
one can form a `subordinated' operator 
$S = \sum_{k\geq 0} F(k) T^k$.    
We obtain asymptotic properties of the subordinated discrete semigroup 
$(S^n\colon  n=1,2, \ldots)$ 
under certain conditions on $F$.  
In particular, we study probabilities $F$ with the property 
that $S$ satisfies the Ritt resolvent condition whenever $T$ 
is power-bounded.    
Examples and counterexamples of this property are discussed.     
The hypothesis of power-boundedness of $T$  
can sometimes be replaced  
by the weaker 
Kreiss resolvent condition.

\bigskip

\bigskip


\section{Introduction}  \label{s1}

Let   
$(S_t)_{t\geq 0}$ 
be a uniformly bounded $C_0$ semigroup 
of operators on 
a complex Banach space $X$.    
Given a family 
$(\mu_t)_{t\geq 0}$ 
of probability measures on   
$\Ri^+: = [0, \infty)$ 
such that 
$\mu_0=\delta_0$ 
(the Dirac measure) and 
$\mu_{t_1} * \mu_{t_2} = \mu_{t_1 + t_2}$,
one can define a new semigroup 
$(\widehat{S}_t)_{t\geq 0}$ 
by 
\begin{equation} 
\widehat{S}_t := \int_0^{\infty} d\mu_t(s) \, S_s 
\label{ecsub} 
\end{equation} 
for 
$t\geq 0$.  
One says that the semigroup  
$(\widehat{S}_t)$ 
is {\it subordinated} to   
$(S_t)$ 
via the convolution semigroup 
$(\mu_t)_{t\geq 0}$.   
For this idea, which goes back to Bochner, see for example  
\cite{CaKa}, \cite[Section 2.4]{Dav80} 
and references therein.

If 
$(\mu_t^{[\alpha]})_{t\geq 0}$ 
is the convolution semigroup of so-called L\'{e}vy 
stable measures of order 
$\alpha \in (0,1)$,      
it is well known that the generator of  
$(\widehat{S}_t)$ is the fractional 
power $-A^{\alpha}$,  
where $-A$ is the generator of the semigroup 
$(S_t=e^{-tA})_{t\geq 0}$.   
In this case a classic result of 
Yosida states that the semigroup 
$(\widehat{S}_t)_{t\geq 0}$ 
is always {\it analytic}; 
see, for example, \cite[Section IX.11]{Yos}.       

Subsequently 
Carasso and Kato \cite{CaKa} asked 
which convolution semigroups 
$(\mu_t)_{t\geq 0}$ on 
$\Ri^+$ 
have the property that the subordinated semigroup is always 
analytic, 
and found various examples and counterexamples of this phenomenon.      

The idea of the present paper is to develop an analogous theory  
for discrete time semigroups.  
Let $X$ be a complex Banach space and 
$T\in {\cal L}(X)$ 
be a bounded linear operator, 
which we usually assume to be {\it power-bounded} in the sense 
that 
$\sup_{n\in \Ni} \| T^n \| < \infty$ 
where 
$\Ni: = \{ 1, 2, 3, \ldots \}$.  
We consider 
$(T^n)_{n\in\Ni}$ as a discrete operator semigroup.    
Let $F$ be a probability on 
$\Zi^+: = \{0, 1, 2, \ldots\}$, 
that is, 
$F\colon \Zi^+ \to [0, 1]$ satisfies 
$\sum_{k\geq 0} F(k)=1$. 
Setting $T^0: = I$,  
we define 
\begin{equation} 
\Psi(F; T): = \sum_{k\geq 0} F(k) T^k \in {\cal L}(X).  
\label{ephidef} 
\end{equation} 
The powers 
$\Psi(F; T)^n$, 
$n\in\Ni$,  
are then given by  
\begin{equation} 
\Psi(F; T)^n = \sum_{k\geq 0} F^{(n)}(k) T^k = \Psi(F^{(n)}; T) 
\label{ediscsub} 
\end{equation} 
where the probability 
$F^{(n)}$    
is the $n$-th convolution power of $F$
(see Section~\ref{s2} below for details).   
Since equation (\ref{ediscsub}) is a discrete analogue of (\ref{ecsub}), 
we think of the discrete semigroup 
$(\Psi(F; T)^n)_{n\in\Ni}$ 
as subordinated to $(T^n)$ 
via the probability $F$.  

We are interested in regularity properties of the subordinated discrete 
semigroup 
$(\Psi(F; T)^n)_{n\in\Ni}$.   
For example, in Section~\ref{s2.5} 
we establish the time regularity property    
\[
\sup_{n\in\Ni} n^{1/2} \| \Psi(F; T)^n - \Psi(F; T)^{n+1} \|
< \infty 
\]
whenever $T$ is power-bounded and 
$F$ satisfies a simple aperiodicity condition.

We will concentrate, however, on a stronger regularity property.  
Let us say that $S\in {\cal L}(X)$ is a {\it Ritt} operator 
if $S$ is power-bounded and 
\[
\sup_{n\in\Ni} n \| S^n - S^{n+1}  \| < \infty. 
\]
This concept  is a discrete time version  of  the notion of a bounded analytic 
$C_0$ semigroup,  
for it is well known 
(for example, \cite[Section 2.5]{Dav80}) 
that a 
$C_0$ semigroup $(S_t = e^{-tA})_{t\geq 0}$ 
is bounded analytic if and only if 
\[
\sup_{t>0} (\| e^{-tA} \| + t \| A e^{-tA} \|) < \infty.  
\]
A basic aim of the current paper is to study probabilities 
on $\Zi^+$ with the property that $\Psi(F; T)$ 
is a Ritt operator for any power-bounded operator $T$.  
Fundamental examples are the probabilities $A_{\alpha}$ 
in the following theorem, 
which are connected with fractional powers of $I-T$  
(see Section~\ref{s5} for details).

\begin{thm} \label{taalpha}  
Let $\alpha\in (0,1)$.  
There exists a probability $A_{\alpha}$ on $\Zi^+$ 
such that for any complex Banach space 
$X$ and any power-bounded operator 
$T\in {\cal L}(X)$, 
the operator 
$\Psi(A_{\alpha}; T) = I - (I-T)^{\alpha}$ 
is a Ritt operator.  
\end{thm} 

Further interesting examples are the `zeta' probabilities $Z_{\alpha}$, 
$\alpha\in (0,1)$.   
These are defined by 
$Z_{\alpha}(0)=0$ and   
\begin{equation} 
Z_{\alpha}(k) : = \zeta(1+\alpha)^{-1} k^{-1-\alpha} 
\label{ezdef} 
\end{equation} 
for $k\in\Ni$, 
where 
$\zeta(s): = \sum_{n\geq 1} n^{-s}$, $s>1$, 
is the zeta function.  
We will establish the following result.  

\begin{thm} \label{tzeta}  
For $\alpha\in (0,1)$ 
and any power-bounded operator 
$T\in {\cal L}(X)$, 
the operator 
$\Psi(Z_{\alpha}; T)$ is a Ritt operator.  
\end{thm} 

Theorems \ref{taalpha} and \ref{tzeta} 
are special cases in Section \ref{s5}, 
where we give useful sufficient (though not necessary) 
conditions on a probability $F$ so that $\Psi(F; T)$ is a Ritt operator 
whenever $T$ is power-bounded. 
Section~\ref{s3} 
discusses some general properties of probabilities with the latter property. 
For example, such probabilities must have an infinite first 
moment,  and therefore have rather slowly decaying tails. 

One can view the theory of large-time behaviour of random 
walks on $\Zi^+$ as behind many of our results. 
On the other hand, underlying the results of Carasso and Kato \cite{CaKa} is the rather different theory 
for  
{\it small-time} analysis of infinitely divisible convolution semigroups on 
$\Ri^+$. 
Thus some of our results on discrete subordinated semigroups 
have a quite different character 
from the continuous case studied in \cite{CaKa}.

Theorem \ref{taalpha} can be extended, by quite different methods, 
to operators $T$ satisfying a weaker condition than 
power-boundedness.  
Namely, say that $T\in {\cal L}(X)$ is a {\it Kreiss} operator 
if it satisfies the Kreiss resolvent condition 
\[
\| (\lambda I - T)^{-1} \| \leq c(|\lambda|-1)^{-1} 
\]
for all $\lambda\in \Ci$ with $|\lambda|>1$.  
The expansion 
$(\lambda I - T)^{-1} = 
\sum_{n\geq 0} \lambda^{-1-n} T^n$ 
($|\lambda|>1$)    
shows that 
every power-bounded operator is a Kreiss operator, 
but the converse is not true in infinite-dimensional Banach spaces 
(see, for example, \cite{Nev3,MSZ} and references therein).      
In general, the convergence of the series (\ref{ephidef}) 
may be 
problematic when $T$ is a Kreiss operator.  
Nevertheless,  
in Section \ref{s7} 
we obtain the following connections between Kreiss operators, 
power-bounded operators and Ritt operators, 
by a different approach through      
the theory of fractional powers of operators.

\begin{thm} \label{tkritt} 
For $T\in {\cal L}(X)$, the following conditions are equivalent. 

\noindent {\rm (I)}$\;$  
$T$ is a Ritt operator. 

\noindent {\rm (II)}$\;$   
There exists a Kreiss operator $S$ and an $\alpha\in (0,1)$   
such that 
$T = I-(I-S)^{\alpha}$.  

\noindent {\rm (III)}$\;$  
There exists a power-bounded operator $S$ and an $\alpha\in (0,1)$ 
such that 
$T = I - (I-S)^{\alpha}$.

\noindent {\rm (IV)}$\;$  
$T$ is power-bounded, 
and there exists a $\gamma_0>1$ 
such that 
$I-(I-T)^{\gamma}$ 
is a Ritt operator for each $\gamma\in (1,\gamma_0)$. 
\end{thm}

The implication (II)$\Rightarrow$(I) in Theorem~\ref{tkritt} 
partially generalizes Theorem~\ref{taalpha}.  

Remark also that condition (IV) in Theorem~\ref{tkritt} 
gives a perturbation result for 
any Ritt operator $T$ under a type of `fractional power'  perturbation.

J. Zem\'{a}nek asked whether there exist Ritt operators $T$ 
with $T\neq I$ and single-point spectrum $\sigma(T)=\{1\}$, 
a question answered affirmatively by Lyubich \cite{Ly}. 
Our results actually reduce Zem\'{a}nek's question 
to an easier question about Kreiss operators, as follows.  

\begin{cor} \label{csingle} 
Suppose that $T\in {\cal L}(X)$ is a Kreiss operator 
with $T\neq I$ and spectrum $\{1\}$. 
Then for each $\alpha\in (0,1)$, 
the operator $S: = I - (I-T)^{\alpha}$ 
is a Ritt operator 
with $S\neq I$ and spectrum $\{1\}$.  
\end{cor}  

For  example, 
consider the Volterra integral 
operator $V$ acting in $X=L^2(0,1)$, 
so 
$(Vf)(x)=\int_0^x dt\, f(t)$ 
for $f\in L^2(0,1)$. 
It is well known that 
the operator  
$T: = (I+V)^{-1}$ 
has spectrum $\{1\}$ and that 
$\|T\|=1$ 
(see \cite[Problem 150]{Hal}, \cite{MSZ}).  
In particular, $T$ is power-bounded, so Corollary~\ref{csingle} applies 
and yields Ritt operators 
$I-(I-T)^{\alpha}$, 
$\alpha\in (0,1)$,  
each with spectrum $\{1\}$.

Actually, 
Lyubich's solution in \cite{Ly} to Zem\'{a}nek's question, 
which involves certain fractional Volterra-type operators,   
also seems a special case of Corollary~\ref{csingle} 
in view of the fact (see \cite{MSZ})
that the operator 
$T_1: = I-V$ 
is a Kreiss operator in the spaces 
$L^p(0,1)$, 
$1\leq p\leq \infty$.

The following fundamental theorem on Ritt operators 
will be an important tool throughout the paper 
(for details of the theorem and further developments 
on Ritt operators, see 
\cite{Blu2,Blu3,Nev2,Nev3,KMOT,Vi,Ly} and their references).   
For $\theta\in [0, \pi)$ 
define the sectors 
$\Lambda_{\theta}: = \{z\in \Ci\colon z\neq 0, |\mathop{\rm Arg} z| < \theta\}$ 
and 
$\overline{\Lambda}_{\theta}
: = \{0\} \cup \{z\in \Ci\colon |\mathop{\rm Arg} z| \leq \theta\}$.   
Put 
$\Di: = \{z\in \Ci\colon |z|<1 \}$ 
and 
$\overline{\Di}: = \{z\in \Ci\colon |z| \leq 1\}$.
In general $\sigma(T)$ denotes the spectrum of $T$.

\begin{thm} \label{tritt} 
For $T\in {\cal L}(X)$, 
the following three conditions (I)-(III) 
are equivalent. 

\noindent {\rm (I)}$\;$  
$T$ is a Ritt operator.

\noindent {\rm (II)}$\;$  
One has 
$\sigma(T) \subseteq \Di \cup \{1\}$; 
and the semigroup 
$(e^{-t(I-T)})_{t\geq 0}$ 
is bounded analytic,  
that is, 
\[
\sup_{t>0} \left( \| e^{-t(I-T)} \| + t \| (I-T) e^{-t(I-T)} \| \right) 
< \infty. 
\]

\noindent {\rm (III)}$\;$ 
$T$ satisfies the Ritt resolvent condition, that is, 
\[
\| (\lambda I -T)^{-1} \| \leq c |\lambda -1|^{-1}
\]
for all $\lambda\in\Ci$ with $|\lambda|>1$. 

\bigskip 

Moreover, if these conditions are satisfied, then there is a $\theta\in (0, \pi/2)$ 
such that 
\[
\sigma(T)  \subseteq (\Di \cup \{1\})  
\cap \{z\in \Ci\colon 1-z\in \overline{\Lambda}_{\theta} \}.  
\]

\end{thm}

The final statement  of Theorem~\ref{tritt} follows from condition (II),  
since for a bounded analytic semigroup 
$(e^{-tA})_{t\geq 0}$ 
it is standard (\cite[Section 2.5]{Dav80}) 
that 
$\sigma(A)\subseteq \overline{\Lambda}_{\theta}$ 
for some $\theta\in (0,\pi/2)$.

\section{Preliminaries}  \label{s2}

This section establishes definitions  
and preliminary results
for subordinated discrete semigroups.  
We refer to 
\cite{Dou,Fol} 
for standard material on Banach algebras and harmonic analysis, 
and to \cite{Fell2,Gal} for Fourier (and other) transforms  
of probability measures.

Suppose that 
${\cal B}$ is   
a complex Banach algebra with unit element 
$1_{\cal B}$.   
The spectrum of an element $A\in {\cal B}$ 
is written  
$\sigma(A; {\cal B})$;    
in case  
${\cal B}={\cal L}(X)$, the algebra of bounded operators 
in a complex Banach space 
$X$, 
this is abbreviated to $\sigma_X(A)$ or $\sigma(A)$.  
For 
$A\in {\cal B}$ 
the exponential is 
$e^A : = \sum_{k\geq 0} (k!)^{-1} A^k \in {\cal B}$, 
where 
$A^0: = 1_{\cal B}$.
One says that $A$ is {\it power-bounded} if 
$\sup_{n\in\Ni} \| A^n \| < \infty$.   
If $A$ is power-bounded then 
$\sigma(A; {\cal B}) \subseteq \overline{\Di}$.

The Banach algebra  $L^1(\Zi)$  consists of all 
functions $F\colon \Zi\to \Ci$ on the integers $\Zi$ 
such that 
$\| F \|_{L^1(\Zi)} : = \sum_{k\in\Zi} |F(k)| < \infty$;  
the (commutative) convolution product 
$F_1 * F_2$  
of $F_1, F_2\in L^1(\Zi)$   
is given by    
\[
(F_1 * F_2)(m) = \sum_{k\in \Zi} F_1(k) F_2(m-k), 
\;\;\; 
m\in\Zi.  
\]   
Write 
$F^{(n)}: = F* F * \cdots * F$ 
for the  
$n$-th convolution power of $F\in L^1(\Zi)$. 
The unit element of $L^1(\Zi)$ is $\delta_0$, 
where in general $\delta_m\in L^1(\Zi)$ is defined so that 
$\delta_m(k)$ is $1$ or $0$ according as $k=m$ or $k\neq m$, 
$k,m \in \Zi$. 

For $F\in L^1(\Zi)$ define the 
$\Zi$-Fourier transform 
$\widehat{F}\in C([-\pi, \pi])$ 
by 
$\widehat{F}(\xi) = \sum_{k\in \Zi} F(k) e^{-ik\xi}$, 
$\xi\in [-\pi, \pi]$. 
Consider the convolution operator $L(F)\in {\cal L}(L^1(\Zi))$ 
given by 
$L(F)F_1: = F * F_1$, 
$F_1\in L^1(\Zi)$. 
Then $L\colon L^1(\Zi) \to {\cal L}(L^1(\Zi))$ 
is an algebra homomorphism, 
with     
\[
\| L(F) \|_{{\cal L}(L^1(\Zi))}  = \| F \|_{L^1(\Zi)}   
\]
and 
\begin{equation} 
\sigma_{L^1(\Zi)}(L(F)) = \sigma(F; L^1(\Zi)) 
= \{ \widehat{F}(\xi) \colon \xi\in [-\pi, \pi] \}
\label{especteq} 
\end{equation} 
for all 
$F\in L^1(\Zi)$. 
For 
$F\in L^1(\Zi)$ define the support   
$\supp(F):= \{k\in\Zi\colon F(k)\neq 0 \}$. 

Put  
\[
L^1(\Zi^+)
 : =  \{F\in L^1(\Zi)\colon \supp(F) \subseteq \Zi^+ = \{0, 1, 2, \ldots \} \};     
\] 
then 
$L^1(\Zi^+)$ is a closed subalgebra of $L^1(\Zi)$.  
For 
$F\in L^1(\Zi^+)$, 
define  
$\phi_F \colon \overline{\Di}\to \Ci$ 
by 
\[
\phi_F(w): = \sum_{k\geq 0} F(k) w^k, 
\;\;\;  
w\in \overline{\Di}.
\] 
Clearly 
$\phi_F$ is continuous on 
$\overline{\Di}$ and 
analytic in 
$\Di$, 
$\| \phi_F \|_{L^{\infty}(\overline{\Di})} \leq \| F \|_{L^1(\Zi)}$,    
and $\widehat{F}$ gives the boundary values of $\phi_F$ in the sense that   
$\phi_F(e^{-i\xi}) = \widehat{F}(\xi)$ for 
$\xi\in [-\pi, \pi]$. 
Moreover 
$\phi_{F_1 * F_2}  =\phi_{F_1} \phi_{F_2}$
for all 
$F_1, F_2\in L^1(\Zi^+)$.

Let 
$\PPi(\Zi): = 
\{F\in L^1(\Zi)\colon F\geq 0, \, \sum_{k\in\Zi} F(k)=1 \}$ 
be the set of probabilities 
on $\Zi$. 
A probability 
$F\in \PPi(\Zi)$ is said to be {\it adapted} if 
$\supp(F)$ 
generates the additive group $\Zi$. 
As a stronger condition,   
$F\in \PPi(\Zi)$ is said to be {\it aperiodic} if all 
of the translated probabilities  
$\delta_m * F$, $m\in \Zi$, 
are adapted.    
It is a useful remark that 
$F\in \PPi(\Zi)$ is aperiodic whenever 
$\supp(F)$ 
contains two consecutive integers $k, k+1$. 

To each $F\in \PPi(\Zi)$ is associated 
a continuous convolution semigroup 
$(e^{-t(\delta_0 - F)})_{t\geq 0} \subseteq \PPi(\Zi)$ 
of probabilities. 
Observe that 
the formula 
$e^{-t(\delta_0-F)} = e^{-t} \sum_{n\geq 0} (n!)^{-1} t^n F^{(n)}$  
shows that these are indeed probabilities.  
 
Set 
$\PPi(\Zi^+) : = \PPi(\Zi) \cap L^1(\Zi^+)$.    
For   
$F\in \PPi(\Zi^+)$, 
the function 
$\phi_F\colon \overline{\Di}\to \overline{\Di}$ 
is called the `generating function' of $F$ 
in classical probability theory.  

Given a power-bounded operator 
$T\in {\cal L}(X)$ and 
any 
$F\in L^1(\Zi^+)$, 
define 
$\Psi(F; T)\in {\cal L}(X)$ by 
equation (\ref{ephidef}). 
Clearly  $\Psi(F; T)$ 
depends linearly on $F$, 
and 
\begin{equation} 
\| \Psi(F; T) \| \leq c(T) \| F \|_{L^1(\Zi)} 
\label{epsinorm} 
\end{equation} 
where 
$c(T): = \sup_{k\geq 0} \| T^k \|$.
One has   
$\Psi(F_1 * F_2; T) = \Psi(F_1; T) \Psi(F_2; T)$ 
for all 
$F_1, F_2\in L^1(\Zi^+)$.  
This yields equation (\ref{ediscsub}) for any 
$F\in L^1(\Zi^+)$.

For 
$F\in L^1(\Zi^+)$ 
it is natural to write  
`$\Psi(F; T) = \phi_F(T)$',  
but we will not use  
the functional calculus notation 
$\phi_F(T)$ 
systematically.    

For any $T\in {\cal L}(X)$ 
(not necessarily power-bounded) 
and any function 
$\varphi$ 
analytic on a neighborhood of 
$\sigma(T)$, 
the operator 
$\varphi(T)\in {\cal L}(X)$ 
is defined by 
the Dunford functional calculus, 
and the spectral mapping theorem 
(see for example 
\cite[Section VIII.7]{Yos}
states that 
$\sigma(\varphi(T)) = \varphi(\sigma(T))$.   
For power-bounded $T$ we have 
$\sigma(T)\subseteq \overline{\Di}$;  
the next result extends the spectral mapping theorem 
to $\Psi(F; T) = \phi_F(T)$ for any $F\in L^1(\Zi^+)$  
(observe that 
$\phi_F$ 
need not be analytic   
on a neighborhood of 
$\overline{\Di}$). 

\begin{thm} \label{tsmap}  
Let $T\in {\cal L}(X)$ be power-bounded.  
Then one has 
$\sigma(\Psi(F; T)) = \phi_F(\sigma(T)) 
\subseteq \phi_F(\overline{\Di})$ 
for all $F\in L^1(\Zi^+)$.   
\end{thm} 
\proof\
Given 
$F\in L^1(\Zi^+)$, 
choose a sequence 
$(F_n)_{n\in\Ni}\subseteq L^1(\Zi^+)$  
such that  
$\lim_{n\to \infty} \| F_n - F \|_{L^1(\Zi)} = 0$ 
and each 
$F_n$ has finite support.  
Since $\Psi(F_n; T)$ is a polynomial in $T$, 
the above-cited spectral mapping theorem yields 
$\sigma(\Psi(F_n; T)) = \phi_{F_n}(\sigma(T))$ 
for all $n$. 
A standard spectral perturbation result 
(\cite[Theorem IV.3.6]{Kat}) 
states that for any commuting operators $S_1, S_2\in {\cal L}(X)$ one has 
\[
\rho(\sigma(S_1), \sigma(S_2)) \leq \| S_1 - S_2 \|, 
\]
where 
\[
\rho(\Gamma_1, \Gamma_2) := \max \left( 
\sup_{x\in \Gamma_1} \inf_{y\in \Gamma_2} |x-y|,  
\sup_{y\in \Gamma_2} \inf_{x\in \Gamma_1} |x-y| \right) 
\]
is the Hausdorff distance between two compact sets 
$\Gamma_1, \Gamma_2 \subseteq \Ci$. 
Since $\lim_{n\to \infty} \| \Psi(F_n; T) - \Psi(F; T) \| = 0$ 
by (\ref{epsinorm}), we deduce that 
\[
\lim_{n\to \infty} \rho(\phi_{F_n}(\sigma(T)), \sigma(\Psi(F; T))) = 0.  
\]
But also
$\lim_{n\to \infty} \rho(\phi_{F_n}(\sigma(T)), \phi_F(\sigma(T))) = 0$, 
as a consequence of the fact that 
$\lim_{n\to \infty} \| \phi_{F_n} - \phi_F \|_{L^{\infty}(\overline{\Di})} 
=0$.  
By uniqueness of limits for the Hausdorff metric  
we obtain 
$\sigma(\Psi(F; T)) = \phi_F(\sigma(T))$.  
\hfill $\Box$ 

\bigskip

The following basic observation is a consequence of (\ref{ediscsub}) 
and (\ref{epsinorm}), 
and  will sometimes be used without comment. 

\begin{lemma} \label{lsubdd} 
Let $T\in {\cal L}(X)$ be power-bounded. 
If $F\in L^1(\Zi^+)$ 
is power-bounded 
(that is, 
$\sup_{n\in \Ni} \| F^{(n)} \|_{L^1(\Zi)} < \infty$),   
then 
$\Psi(F; T)\in {\cal L}(X)$ 
is power-bounded.   
In particular, $\Psi(F; T)$ is power-bounded for each
$F\in \PPi(\Zi^+)$. 
\end{lemma}

\section{A basic regularity result}  \label{s2.5}

In this section we establish a basic regularity estimate 
for the differences 
$\Psi(F; T)^n - \Psi(F; T)^{n+1}$. 

We first quote a result which characterizes a decay of order  
$O(n^{-1/2})$ of the 
differences 
$\| S^n-S^{n+1}\|$,  
for an 
$S\in {\cal L}(X)$; 
note that 
this decay is slower than for Ritt operators.

\begin{thm} \label{tnevhalf} 
Given $S\in {\cal L}(X)$, the following two conditions are equivalent. 

\noindent {\rm (I)}$\;$  
The operator $S$ is power-bounded, 
and 
$\sup_{n\in\Ni} n^{1/2} \| S^n - S^{n+1} \| < \infty$.

\noindent {\rm (II)}$\;$ 
There exist $\beta\in (0,1)$ and a power-bounded operator 
$T\in {\cal L}(X)$ such that 
$S = \beta T + (1-\beta) I$. 

\bigskip 

If these conditions hold then there is a $\beta\in (0,1)$ 
such that 
\[
\sigma(S)\subseteq \{z\in \Ci\colon |z-(1-\beta)| \leq \beta \} 
\subseteq \Di\cup \{1\}.      
\] 
\end{thm} 

The implication (II)$\Rightarrow$(I) in Theorem~\ref{tnevhalf} 
is  contained in \cite[Lemma~2.1]{FW}, 
\cite[Theorem~4.5.3]{Nev1} 
or \cite[Theorem~8]{Nev3}. 
The implication (I)$\Rightarrow$(II) was recently proved by the author 
(see \cite{Dun10}, where a number of other equivalent conditions on $S$ 
are studied).  
However, strictly speaking we do not need the latter implication 
in the present paper. 

The statement about 
$\sigma(S)$ in Theorem~\ref{tnevhalf} is a consequence 
of condition (II) and the fact that 
$\sigma(T)\subseteq \overline{\Di}$ 
for power-bounded $T$.

In Theorem~\ref{tnevhalf}(II) we have  
$S= \Psi(F; T)$,     
where 
$F\in \PPi(\Zi^+)$ is the `Bernoulli' 
probability with 
$\supp(F)=\{ 0,1 \}$,   
$F(0)=1-\beta$, 
$F(1)=\beta$.  
The main result  of this section, stated next,   
generalizes the implication (II)$\Rightarrow$(I) 
in Theorem~\ref{tnevhalf} to  
arbitrary aperiodic probabilities.

\begin{thm} \label{taphalf}  
Let $T\in {\cal L}(X)$ be power-bounded, 
and let 
$F\in \PPi(\Zi^+)$ 
be aperiodic. 
Then
\begin{equation} 
\sup_{n\in\Ni} n^{1/2} \|  \Psi(F; T)^n - \Psi(F; T)^{n+1} \| 
< \infty, 
\label{esubhalf} 
\end{equation} 
and there exists a
$\beta\in (0,1)$ 
such that  
\[
\sigma(\Psi(F; T)) 
\subseteq \{z\in \Ci\colon |z-(1-\beta)| \leq \beta \} 
\subseteq \Di \cup \{1\}. 
\]  
\end{thm}

To prove Theorem~\ref{taphalf} we require the following
result proved 
in the Appendix 
(Section~\ref{sapp} below).  
The equivalence of conditions (I) and (II) below  
is known 
(see  \cite[Chapter II]{Spi}, 
where the term `strongly aperiodic' is used to mean aperiodic).   
However, to our knowledge,  condition (IV) seems to be new and interesting.

\begin{prop} \label{papchar} 
Let 
$F\in \PPi(\Zi)$ be adapted. 
Then the following conditions are equivalent. 

\noindent {\rm (I)}$\;$  
$F$ is aperiodic. 

\noindent {\rm (II)}$\;$  
$| \widehat{F}(\xi)|< 1$ 
for all $\xi \in [-\pi, \pi]\backslash \{0\}$. 

\noindent {\rm (III)} $\;$   
$\widehat{F}(\xi) \in  \Di \cup \{1\}$ 
for all 
$\xi\in [-\pi, \pi]$.

\noindent {\rm (IV)}$\;$  
There exist 
$\beta\in (0,1)$ and a power-bounded element  
$G\in L^1(\Zi)$ 
(that is, 
$\sup_{n\in\Ni} \| G^{(n)} \|_{L^1(\Zi)} < \infty$)  
such that 
$F = \beta G + (1-\beta) \delta_0$. 

\bigskip

If these conditions hold, 
then 
\begin{equation} 
\sup_{n\in \Ni} n^{1/2} \| F^{(n)} - F^{(n+1)} \|_{L^1(\Zi)} < \infty. 
\label{eaphalf} 
\end{equation}
\end{prop}

\noindent {\bf Proof of Theorem~\ref{taphalf}.}  
Choose $\beta$ and $G$ 
as in 
part (IV) of Proposition~\ref{papchar}.  
Then 
$G\in L^1(\Zi^+)$,  
the operator $\Psi(G; T)$ 
is power-bounded, 
and      
$\Psi(F; T) = \beta \Psi(G; T) + (1-\beta) I$.
Theorem~\ref{taphalf} therefore follows from Theorem~\ref{tnevhalf}. 

An alternative way to obtain (\ref{esubhalf}) is to write   
\[
\Psi(F; T)^n - \Psi(F;  T)^{n+1} = 
\Psi( F^{(n)} - F^{(n+1)}; T)
\] 
and apply (\ref{eaphalf}) together with  
(\ref{epsinorm}).  
\hfill $\Box$

\bigskip

\noindent {\bf Remark.  }  
If
$F\in \PPi(\Zi)$ 
satisfies 
$F(0)>0$ 
then condition (IV) of Proposition~\ref{papchar} is easy to obtain 
(and aperiodicity is not needed).     
In fact, it suffices to choose $\beta\in (1-F(0), 1)$ 
and observe that 
\[
G: = \beta^{-1} (F - (1-\beta)\delta_0) 
\] 
is an element of $\PPi(\Zi)$. 
However, 
if $F(0)=0$ 
(a case we shall often encounter),   
then $G$ in condition (IV) cannot be a probability.

\section{The class ${\cal A}$} \label{s3}

In this section we begin our study of 
the class  
${\cal A}= {\cal A}(\Zi^+)  \subseteq \PPi(\Zi^+)$ 
of `Ritt' probabilities,  
which we define by   
\[
{\cal A}: = \{ F\in \PPi(\Zi^+)\colon 
 \sup_{n\in\Ni} n \| F^{(n)} - F^{(n+1)} \|_{L^1(\Zi)}  
< \infty \}. 
\]
Clearly 
$\delta_0 \in {\cal A}$;   
non-trivial examples of probabilities in ${\cal A}$ will be constructed 
in Sections~\ref{s5} and \ref{s6}.           
In this section we develop general properties of the class ${\cal A}$, 
including its connection with 
subordinated discrete semigroups.  

Here is our basic characterization of the class
${\cal A}$.  
It was partly inspired by 
the results on continuous  
convolution semigroups in 
\cite[Theorem 4]{CaKa}.

\begin{thm} \label{tabaschar} 
Given 
$F\in \PPi(\Zi^+)$, 
the following four conditions are equivalent. 

\noindent {\rm (I)}$\;$  
$F\in {\cal A}$. 

\noindent {\rm (II)}$\;$  
The operator $L(F)$ is a Ritt operator in 
$L^1(\Zi)$. 

\noindent {\rm (III)}$\;$ 
For any complex Banach space 
$X$ and any power-bounded operator 
$T\in {\cal L}(X)$, 
the operator $\Psi(F; T)$ is a Ritt operator.    

\noindent {\rm (IV)}$\;$  
One has $\widehat{F}(\xi) \in \Di\cup \{1\}$ 
for all $\xi\in [-\pi, \pi]$, 
and  
\[
\sup_{t\geq 1} t \| (\delta_0 -F) * e^{-t(\delta_0 -F)} \|_{L^1(\Zi)}  
< \infty.    
\]

\bigskip 

Moreover, if 
$F\in {\cal A}$,  
then there is a 
$\theta\in (0,\pi/2)$  
such that 
\begin{equation} 
\phi_F(w) \in \Di \cup \{1\}, 
\;\;\; 
1-\phi_F(w) \in \overline{\Lambda}_{\theta} 
\label{esectf} 
\end{equation} 
for all 
$w\in \overline{\Di}$.  
Finally, if 
$F\in {\cal A}\backslash \{ \delta_0 \}$ 
then  
\begin{equation} 
\sum_{k\geq 0} k F(k) = + \infty  
\label{einfmean} 
\end{equation}  
and hence $\supp(F)$ is infinite.  
\end{thm} 
\proof\
(I)$\Leftrightarrow$(II): 
observe that 
\[
\| L(F)^n - L(F)^{n+1} \|_{{\cal L}(L^1(\Zi))}  
 = \| F^{(n)} - F^{(n+1)} \|_{L^1(\Zi)}. 
\]

(I)$\Rightarrow$(III):     
setting
$c(T): = \sup_{k\geq 0} \| T^k \|$, 
one has  
\[
\| \Psi(F;T)^n - \Psi(F; T)^{n+1} \| 
= \| \Psi(F^{(n)} - F^{(n+1)}; T) \| 
\leq c(T)  \| F^{(n)} - F^{(n+1)} \|_{L^1(\Zi)} 
\]
for all $n\in\Ni$, 
where the last step used   
(\ref{epsinorm}).   

(III)$\Rightarrow$(II): 
consider the operator  
$T: = L(\delta_1)$ 
acting in 
$X=L^1(\Zi)$, 
and  observe that  
\[
L(F) = \sum_{k\geq 0} F(k) T^k = \Psi(F; T). 
\]

(IV)$\Leftrightarrow$(II): 
recalling (\ref{especteq}), 
one sees that 
condition (IV) holds if and only if  
$L(F)\in {\cal L}(L^1(\Zi))$ 
satisfies the conditions   
$\sigma(L(F)) \subseteq \Di\cup \{1\}$  
and 
\[
\sup_{t\geq 1} 
t \| (I-L(F)) e^{-t(I-L(F))} \| 
 < \infty.  
\]   
By Theorem~\ref{tritt}, 
those conditions hold if and only if 
$L(F)$ is a Ritt operator in  
$L^1(\Zi)$.  

Thus we have shown the equivalence of all the conditions (I)-(IV).

We next assume that 
$F\in {\cal A}$ 
and establish (\ref{esectf}).  
Consider the Banach space 
$X=C(\overline{\Di})$ of continuous functions on $\overline{\Di}$ 
under the usual supremum norm,  
and let 
$T\in {\cal L}(X)$ 
be the multiplication operator 
$(Tf)(w): = wf(w)$ 
for 
$f\in C(\overline{\Di})$, 
$w\in \overline{\Di}$. 
Then 
$\|T \| =1$,  
and $\Psi(F; T)$ is the multiplication operator   
\[
(\Psi(F; T)f)(w) = \phi_F(w) f(w), 
\;\;\; 
 w\in \overline{\Di}.  
\]
By applying the  final statement of Theorem~\ref{tritt} to the Ritt  
operator
$\Psi(F; T)$, 
since 
$\sigma(\Psi(F; T)) = \phi_F(\overline{\Di})$      
we obtain (\ref{esectf}).

To prove the last assertion (\ref{einfmean}), 
suppose that 
$F\in \PPi(\Zi^+)$, 
$F \neq \delta_0$,  
and 
$a: = \sum_{k\geq 0} k F(k) < \infty$; 
we will show that 
$F \notin {\cal A}$.  
Note first that $a>0$. 
Differentiating 
$\widehat{F}(\xi) = \sum_{k\geq 0} F(k) e^{-ik\xi}$ 
with respect to $\xi$ 
shows that 
$\widehat{F} \in C^1([-\pi, \pi])$ 
and 
$\widehat{F}'(0) = -ia$. 
Since $\widehat{F}(0)=1$ 
we have 
\[
\widehat{F}(\xi) = 1 - ia\xi + o(|\xi|) 
\]
where 
$o(|\xi|)$ denotes   
a function of 
$\xi \in [-\pi, \pi]$ 
such that 
$\lim_{\xi\to 0} o(|\xi|)/|\xi| = 0$. 
Since $a>0$, 
for some small $\varepsilon\in (0,1)$ 
we have 
$| 1 - \widehat{F}(\xi) | \geq \varepsilon |\xi|$ 
for all 
$\xi\in [-\varepsilon, \varepsilon]$, 
while  
$1-\RRe \widehat{F}(\xi) = o(|\xi|)$.
Because   
$\widehat{F}(\xi) = \phi_F(e^{-i\xi})$, 
it follows that the second statement of 
(\ref{esectf}) cannot hold.    
Therefore 
$F\notin {\cal A}$.    
\hfill $\Box$ 

\bigskip

\noindent {\bf Remarks.} 
For the Laplace transforms of certain continuous convolution semigroups 
on $\Ri^+$, 
\cite[Theorem~4]{CaKa} gives 
a sectorial condition somewhat  analogous to the second part of (\ref{esectf}).

Although we believe that the conditions (\ref{esectf}) are not {\it sufficient} 
for a probability $F$ to be in ${\cal A}$, 
it does not seem easy to find counterexamples.

Regarding (\ref{einfmean}),  
for   
$F\in \PPi(\Zi)$  
the quantity 
$\sum_{k\in \Zi} k F(k)$ 
(when it exists) 
is the first moment or `center of mass' of 
$F$.     
On any compactly generated locally compact group 
(in particular on $\Zi$) 
there is a large class of `Ritt' probability measures which have
finite, in fact vanishing, first moments;   
see \cite{Dun3} for this theory.  
On $\Zi^+$,  however, 
the only probability with a vanishing first moment is $\delta_0$.  

Condition (\ref{einfmean}) shows that many commonly studied probabilities on $\Zi^+$ 
(for example, probabilities of finite support, 
or the Poisson probability 
$P_s: = e^{-s} \sum_{k\geq 0} (k!)^{-1} s^k \delta_k$ 
for any 
$s>0$)     
are not elements of 
${\cal A}$.     

In Section~\ref{s6} we shall describe an example satisfying (\ref{einfmean}) 
but not (\ref{esectf}).

\bigskip

Here are some further interesting properties of the class ${\cal A}$. 

\begin{prop} \label{paprop} 
The set 
${\cal A}$ 
is a convex subset of 
$\PPi(\Zi^+)$. 
For any $F_1, F_2\in {\cal A}$ 
and $G\in \PPi(\Zi^+)$, 
one has 
$F_1 * F_2 \in {\cal A}$ and 
\begin{equation} 
\sum_{k\geq 0} F_1(k) G^{(k)} \in {\cal A}. 
\label{easubord} 
\end{equation} 
\end{prop} 
\proof\
Let $F_1, F_2\in {\cal A}$. 
From 
\[
(F_1*F_2)^{(n)} - (F_1*F_2)^{(n+1)} 
= (F_1^{(n)} - F_1^{(n+1)}) * F_2^{(n)} 
  + F_1^{(n+1)} * ( F_2^{(n)} - F_2^{(n+1)}) 
\]
it easily follows that $F_1*F_2 \in {\cal A}$. 
For a convex combination 
$F= \lambda F_1+(1-\lambda) F_2 \in \PPi(\Zi^+)$  
where 
$\lambda \in (0,1)$, 
one has  
\[
e^{-t(\delta_0-F)} = e^{-t\lambda(\delta_0- F_1)} * e^{-t(1-\lambda)(\delta_0 - F_2)}.  
\]
By differentiating this formula with respect to 
$t$ and applying Condition (IV) of 
Theorem~\ref{tabaschar}, 
one sees that 
$F\in {\cal A}$. 
Thus   
${\cal A}$ is convex.  
Finally, 
consider the probability   
$H: = \sum_{k\geq 0} F_1(k) G^{(k)} \in \PPi(\Zi^+)$.  
Since 
$L(H)= \Psi(F_1; L(G))$ 
and 
$F_1\in {\cal A}$,  
we deduce from Theorem~\ref{tabaschar}(III) that 
$L(H)$ is a Ritt operator in  
$L^1(\Zi)$, 
hence $H\in {\cal A}$. 
\hfill $\Box$

\bigskip

It follows from 
Theorem~\ref{tabaschar}(IV) and Proposition~\ref{papchar} that 
an element of ${\cal A}$ is aperiodic if and only if it is adapted. 
The next result shows that non-aperiodic elements of ${\cal A}$ 
can be obtained by scaling of aperiodic elements.   
Thus in practice it suffices to study the aperiodic elements of 
${\cal A}$.    

\begin{lemma} \label{lapred} 
Suppose that $F\in {\cal A}\backslash \{\delta_0 \}$ 
is not aperiodic. 
Then there exist 
$m\in \{2, 3, 4, \ldots \}$ 
and an aperiodic probability 
$\widetilde{F} \in {\cal A}$ 
such that 
$\supp(F) \subseteq \{0, m, 2m, \ldots \}$ 
and 
$F(km) = \widetilde{F}(k)$ 
for all $k\in \Zi^+$. 
\end{lemma} 
\proof\
By the remark preceding the lemma, 
$F$ is not adapted.  
Since 
$F\neq \delta_0$, 
the subgroup of $\Zi$ generated by 
$\supp(F)$ must equal  
$m\Zi$ for some 
$m\in \{2, 3, \ldots \}$. 
Set 
$\widetilde{F}(k): = F(km)$ for 
$k\in \Zi^+$. 
Then 
$\widetilde{F}\in \PPi(\Zi^+)$ is adapted, 
and 
$\widetilde{F} \in {\cal A}$ 
because 
$\| \widetilde{F}^{(n)} - \widetilde{F}^{(n+1)} \|_{L^1(\Zi)} 
= \| F^{(n)} - F^{(n+1)} \|_{L^1(\Zi)}$ 
for all 
$n\in\Ni$. 
By the remark preceding the lemma, 
$\widetilde{F}$ 
is aperiodic.  
\hfill $\Box$

\section{Sufficient conditions and examples} \label{s5}

In this section 
we first give a sufficient criterion 
for a probability to belong to 
${\cal A}$, 
in terms of the $\Zi$-Fourier 
transform. 
This result leads to fundamental examples of elements of ${\cal A}$, 
including the probabilities $A_{\alpha}$, 
$Z_{\alpha}$ 
in Theorems \ref{taalpha} and \ref{tzeta}.   
 
Write  
$\partial_{\xi}$ 
for differentiation with respect to the real variable  
$\xi$.

\begin{thm} \label{tcharsuff} 
Let $\alpha\in (0,1)$.  

\noindent {\rm (I)}$\;$  
Suppose that 
$F\in \PPi(\Zi^+)$ is aperiodic 
and there exist 
$\varepsilon\in (0,1)$, 
$c>0$ such that    
\begin{equation} 
1 - \RRe \widehat{F}(\xi) 
\geq \varepsilon  |\xi|^{\alpha} 
\label{ereal} 
\end{equation}
for all 
$\xi\in [-\varepsilon, \varepsilon]$ 
and 
\begin{equation} 
| \partial_{\xi} \widehat{F}(\xi) | \leq c |\xi|^{\alpha-1} 
\label{efderiv} 
\end{equation} 
for all $\xi\in [-\pi, \pi]\backslash \{0\}$. 
Then $F\in {\cal A}$.

\bigskip 

\noindent {\rm (II)}$\;$  
Let $F, G\in \PPi(\Zi^+)$ 
be aperiodic  
with $F$ satisfying the conditions 
of part {\rm (I)} above.  
If there exists an  
$a>0$ such that 
\[
\sum_{k\geq 1} k |G(k) - a F(k)| < \infty, 
\]
then $G$ satisfies the conditions of part {\rm (I)} above, 
and hence  
$G\in {\cal A}$.   
\end{thm} 

Part (I) above is partly  
inspired by an  argument of Carasso and Kato 
dealing with Fourier transforms of L\'{e}vy stable measures    
(see \cite[Example 3]{CaKa}).

\bigskip

\noindent {\bf Proof of Theorem~\ref{tcharsuff}.}   
To prove part (I), we will verify
condition (IV) of Theorem~\ref{tabaschar}. 
By Proposition~\ref{papchar} one has 
$\widehat{F}(\xi) \subseteq \Di\cup \{1\}$ 
for all 
$\xi\in [-\pi, \pi]$. 
Moreover, 
$|\widehat{F}(\xi)|<1$ when  
$\xi\neq 0$, 
so by modifying   
$\varepsilon$ 
we may assume that     
(\ref{ereal}) 
holds for all 
$\xi\in [-\pi, \pi]$.

Set 
$G_t: =  (\delta_0 - F) * e^{-t(\delta_0-F)} 
\in L^1(\Zi^+)$ 
for each 
$t>0$.  
Put 
$R_t: = \widehat{G_t}$, so that 
\[
R_t(\xi) = (1-\widehat{F}(\xi)) e^{-t(1-\widehat{F}(\xi))}. 
\]
Integrating (\ref{efderiv}) 
yields a bound   
$|1- \widehat{F}(\xi)| = |\widehat{F}(0) - \widehat{F}(\xi)|    
 \leq c' |\xi|^{\alpha}$
for 
$\xi\in [-\pi, \pi]$. 
Using this together with (\ref{ereal}), 
one finds that there are 
$c,b>0$ 
such that  
\begin{eqnarray} 
|R_t(\xi)| 
& = & |1-\widehat{F}(\xi)| e^{-t(1-\RRe \widehat{F}(\xi))} 
\leq c |\xi|^{\alpha} e^{-b t |\xi|^{\alpha}},  
 \nonumber \\ 
|\partial_{\xi} R_t(\xi)| 
& \leq & c |\xi|^{\alpha-1} e^{-b t |\xi|^{\alpha}}
 \label{ertbd}  
\end{eqnarray} 
for all 
$t\geq 1$ and  
$\xi\in [-\pi, \pi]$. 
Note that 
$\partial_{\xi} R_t$ 
is the $\Zi$-Fourier transform 
of the function 
$H_t(k): = -ik G_t(k)$, 
$k\in \Zi$.  
 
The next part of the argument is similar to   
\cite[pp.875-876]{CaKa}. 
Fix $r>1$ with  
$r< \min\{2, (1-\alpha)^{-1} \}$, 
and set 
$s:=r/(r-1) \in (2, \infty)$. 
Let 
$\| \cdot \|_r$ 
denote the norm in the space  
$L^r([-\pi, \pi]; d\xi)$.    
Integrating (\ref{ertbd}) gives bounds of type  
\[
\| R_t \|_r \leq c t^{-1} t^{-1/(\alpha r)}, 
\;\;\; 
\| \partial_{\xi} R_t \|_r \leq c t^{-1} t^{1/(\alpha s)} 
\]
for all 
$t\geq 1$. 
The inverse $\Zi$-Fourier transform is a bounded operator  
from $L^r([-\pi, \pi])$ 
into $L^s(\Zi)$, 
by the Hausdorff-Young inequality 
(\cite[Theorem 12.12]{Rudrca}); 
hence 
\[
\| G_t \|_{L^s(\Zi)}  \leq c' t^{-1} t^{-1/(\alpha r)}, 
\;\;\; 
\| H_t  \|_{L^s(\Zi)}  \leq c' t^{-1} t^{1/(\alpha s)} 
\]
for $t\geq 1$. 
Using H\"{o}lder's inequality we obtain 
\begin{eqnarray*} 
\| G_t \|_{L^1(\Zi)} 
& \leq & 
  \left( \sum_{k\geq 0} (1+ kt^{-1/\alpha})^{-r} \right)^{1/r}   \\
 & & {} \times 
 \left( \sum_{k\geq 0} |G_t(k)|^s (1+kt^{-1/\alpha})^s \right)^{1/s} \\ 
 & \leq & c t^{1/(\alpha r)} \left[ \| G_t \|_{L^s(\Zi)} 
+ t^{-1/\alpha} \| H_t \|_{L^s(\Zi)} \right] 
 \leq c' t^{-1} 
\end{eqnarray*} 
for all $t\geq 1$. 
Thus  
$F$ satisfies Condition (IV) of Theorem~\ref{tabaschar}, 
and  
$F\in {\cal A}$.

To prove part (II), 
put  
$P: = G - a F \in L^1(\Zi^+)$. 
The hypothesis 
$\sum_{k\geq 1} k |P(k)| < \infty$ 
implies that 
$\widehat{P} 
\in C^1([-\pi, \pi])$.  
Since 
$\widehat{G}(0) = \widehat{F}(0) =1$
one has 
$\widehat{P}(0)= 1-a$,  
and hence 
$| \widehat{P}(\xi) - (1-a) | \leq c' |\xi|$ 
for all 
$\xi\in [-\pi, \pi]$. 
The equations 
\[
1-\widehat{G} = a (1 -\widehat{F}) - (\widehat{P}-(1-a)), 
\;\;\; 
\partial_{\xi} \widehat{G}(\xi)  
= a \partial_{\xi} \widehat{F}(\xi)  + \partial_{\xi}  \widehat{P}(\xi) 
\]
then show that $G$ satisfies estimates of type  
(\ref{ereal}) and  (\ref{efderiv}).      
\hfill $\Box$ 

\bigskip

Our first application of Theorem~\ref{tcharsuff}, 
to probabilities 
$A_{\alpha}$, $\alpha\in (0,1)$,   
is based on the following lemma.    
Write 
\[
\left( \begin{array}{c} \beta \\ k \end{array} \right) 
: = (k!)^{-1} \beta (\beta-1) \cdots (\beta-(k-1))
\]
for 
$\beta\in \Ri$ 
and $k\in \Ni$. 
For $\beta>0$ we always consider 
the principal branch 
of the power function $z\mapsto z^{\beta}$ 
defined for $z\in \Ci\backslash (-\infty, 0)$.

\begin{lemma} \label{lalem} 
Given $\alpha\in (0,1)$,  
set  
\begin{equation} 
a_k = a_k^{(\alpha)} 
:= (-1)^{k-1} 
\left( \begin{array}{c} \alpha \\ k \end{array} \right) 
= \frac{\alpha (1-\alpha) (2-\alpha) \cdots ((k-1)-\alpha)}{k!}   
\label{eadef} 
\end{equation}
for 
$k\in\Ni$. 
Then
$a_k>0$ for all 
$k\in\Ni$, 
$\sum_{k\geq 1} a_k =1$, 
and  
\begin{equation} 
1 - (1-w)^{\alpha} = \sum_{k\geq 1} a_k w^k 
\label{efracgen} 
\end{equation} 
for all 
$w\in \overline{\Di}$.  
In addition,  
\begin{equation} 
a_k = \alpha \Gamma(1-\alpha)^{-1} k^{-\alpha-1} 
 + O(k^{-\alpha-2}) 
\label{efasymp} 
\end{equation} 
for 
$k\in\Ni$. 
\end{lemma} 
\proof\
The binomial series expansion   
\[
(1-w)^{\alpha} = 1 + \sum_{k\geq 1} 
\left( \begin{array}{c} \alpha \\ k \end{array} \right) 
(-1)^k w^k 
\]
is valid for all 
$w\in \Di$, 
giving (\ref{efracgen}) 
for 
$w\in \Di$.    
The last expression in (\ref{eadef}) shows that 
$a_k>0$. 
For any $r\in (0,1)$ and 
$N\in\Ni$ 
we have  
\[
\sum_{k=1}^N a_k r^k \leq \sum_{k\geq 1} a_k r^k = 1 -(1-r)^{\alpha} < 1.  
\]
Taking limits as $r\to 1$ and then as 
$N\to \infty$ yields   
$\sum_{k\geq 1} a_k \leq 1$.      
Thus each side of 
(\ref{efracgen}) defines a continuous function of 
$w\in \overline{\Di}$, 
so by continuity (\ref{efracgen}) holds for all 
$w\in {\overline{\Di}}$.  
Taking 
$w=1$ in (\ref{efracgen}) 
yields 
$\sum_{k\geq 1} a_k=1$.

As a preliminary step to (\ref{efasymp}), we claim that 
\begin{equation} 
a_k 
\sim \alpha \Gamma(1-\alpha)^{-1} k^{-\alpha-1}, 
\;\;\; k\to \infty, 
\label{efasymp0} 
\end{equation} 
where 
the notation 
$\sim$ means 
that the ratio of the two sides converges to 
$1$ as 
$k\to \infty$. 
Now
\[
a_k 
= \frac{\alpha ((k-1)-\alpha)((k-2)-\alpha) \cdots (1-\alpha)
\Gamma(1-\alpha)}{k! \Gamma(1-\alpha)}   
= \frac{\alpha \Gamma(k-\alpha)}{\Gamma(1-\alpha) \Gamma(k+1)}   
\]
and from Stirling's formula
$\Gamma(t+1) \sim (2\pi t)^{1/2} t^t e^{-t}$, $t\to \infty$,  
one has    
\begin{eqnarray*} 
\frac{\Gamma(k-\alpha)}{\Gamma(k+1)} 
& \sim & \frac{(k-\alpha-1)^{k-\alpha-1 + (1/2)} e^{-(k-\alpha-1)}}
         {k^{k+ (1/2)} e^{-k}}    \\
& = & e^{\alpha+1}(1-(\alpha+1)k^{-1})^{k+(1/2)} (k-\alpha-1)^{-\alpha-1} 
       \\
& \sim & (k-\alpha-1)^{-\alpha-1} \sim k^{-\alpha-1}     
\end{eqnarray*}  
from which (\ref{efasymp0}) follows. 
Using the precise 
Stirling's estimate  
\[
\left| \frac{\Gamma(t+1)}{(2\pi t)^{1/2} t^t e^{-t}} - 1 \right| 
\leq c t^{-1}, 
\;\;\; 
t\geq 1,  
\] 
it is straightforward to refine the proof of (\ref{efasymp0}) 
to get (\ref{efasymp}); 
we leave this to the reader. 
\hfill $\Box$

\bigskip

Applying Lemma~\ref{lalem}, 
for each $\alpha\in (0,1)$  
define  
$A_{\alpha}\in \PPi(\Zi^+)$ 
by 
$A_{\alpha}(0): = 0$ and 
\begin{equation} 
A_{\alpha}(k) : = a_k^{(\alpha)} 
= \frac{\alpha (1-\alpha) (2-\alpha) \cdots ((k-1)-\alpha)}{k!} 
\label{efdens}  
\end{equation} 
for $k\in\Ni$. 
Lemma~\ref{lalem} implies that $A_{\alpha}$ 
has generating function  
\begin{equation} 
\phi_{A_{\alpha}}(w) = 1 - (1-w)^{\alpha}. 
\label{eagen} 
\end{equation} 
Given a power-bounded operator $T$, 
it is therefore natural to write 
$\Psi(A_{\alpha}; T) = I - (I-T)^{\alpha}$.
Since  
\[
\widehat{A}_{\alpha}(\xi) = 1 - (1- e^{-i\xi})^{\alpha} 
\]
for $\xi\in [-\pi, \pi]$,   
one easily checks that 
$F=A_{\alpha}$ 
satisfies the conditions of Theorem~\ref{tcharsuff}(I),    
so that 
$A_{\alpha}\in {\cal A}$. 
These observations and Theorem~\ref{tabaschar} 
yield Theorem~\ref{taalpha}, 
which we now restate more precisely. 
 
\begin{thm} \label{texalpha} 
For each $\alpha\in (0,1)$,    
the probability 
$A_{\alpha}\in \PPi(\Zi^+)$ 
defined by (\ref{efdens})     
is an element of ${\cal A}$.  
For any power-bounded operator $T\in {\cal L}(X)$, 
the operator $\Psi(A_{\alpha}; T) = I-(I-T)^{\alpha}$ 
is a Ritt operator.      
\end{thm}

We next consider the zeta probabilities $Z_{\alpha}\in \PPi(\Zi^+)$ 
defined by (\ref{ezdef}).  

\begin{thm} \label{tzeta2} 
For each $\alpha\in (0,1)$ one has $Z_{\alpha}\in {\cal A}$. 
\end{thm} 

Observe that Theorem~\ref{tzeta2} yields Theorem~\ref{tzeta}. 

\bigskip 

\noindent {\bf Proof of Theorem~\ref{tzeta2}.}      
We mention two methods of proof.       
The first method is to observe from   
(\ref{efasymp})
that there is a constant $c_{\alpha}>0$ with  
\[
Z_{\alpha}(k) = c_{\alpha} A_{\alpha}(k) + O (k^{-2-\alpha}), 
\;\;\; 
k\in \Ni.   
\]
Applying Theorem~\ref{tcharsuff}(II) with $F=A_{\alpha}$, 
$G=Z_{\alpha}$  
shows that  
$Z_{\alpha} \in {\cal A}$.  

The second method involves the polylogarithm 
$\mathop{\rm Li}_s$, 
which is  
defined by 
\[
{\mathop{\rm Li}}_s(w): = \sum_{k\geq 1} k^{-s} w^k 
\]
for  
$s>1$, 
$w\in \overline{\Di}$,  
and  can be analytically continued in 
$w$ to a function analytic 
on 
$\Ci\backslash [1,\infty)$.   
For 
$s>1$ with 
$s\notin \Ni$,  
one has a series expansion 
\begin{equation} 
{\mathop{\rm Li}}_s(e^{\mu}) = \zeta(s) + \Gamma(1-s) (-\mu)^{s-1} 
+ \sum_{n\geq 1} \zeta(s-n) (n!)^{-1} \mu^n
\label{eplexp} 
\end{equation} 
for  all 
$\mu\in \Ci$ with 
$|\mu|< 2\pi$,  
$\mu \notin (0,2\pi)$.    
This result is contained in the expansion of the special function 
$\Phi$ 
in \cite[p.29,equation (8)]{EMOT}, 
noting that 
$\mathop{\rm Li_s}(w)= w \Phi(w,s,1)$.   
Taking 
$s=1+\alpha$ and   
$\mu=-i\xi$ in (\ref{eplexp}) 
one obtains 
\begin{eqnarray*} 
\widehat{Z}_{\alpha}(\xi) 
& = & \zeta(1+\alpha)^{-1} 
{\mathop{\rm Li}}_{1+\alpha}(e^{-i\xi})   \\ 
& = & 1 - b_{\alpha} (i\xi)^{\alpha} + \xi \varphi_{\alpha}(\xi)
\end{eqnarray*}   
for all 
$\xi\in [-\pi, \pi]$, 
where 
$b_{\alpha}: = -\Gamma(-\alpha) \zeta(1+\alpha)^{-1} >0$
and where 
$\varphi_{\alpha}\in C^{\infty}([-\pi, \pi])$ 
is a smooth function.  
It readily follows that 
$F=Z_{\alpha}$ 
satisfies 
the conditions of Theorem~\ref{tcharsuff}(I), 
so $Z_{\alpha}\in {\cal A}$.  
\hfill $\Box$

\bigskip 

The next result considerably generalizes the examples 
$Z_{\alpha}$ 
and $A_{\alpha}$.

\begin{thm} \label{tpexamples}  
Suppose that 
$F\in \PPi(\Zi^+)$ is such that   
\[ 
F(k) = \sum_{j=1}^N c_j k^{-1-\alpha_j} + P(k) 
\]
for all $k\in\Ni$,   
where 
$N\in\Ni$, 
$0< \alpha_1 < \cdots < \alpha_N < 1$,  
$c_1, \ldots, c_N>0$, 
and 
$P\in L^1(\Zi^+)$ 
with  
$\sum_{k\geq 1} k |P(k)| < \infty$. 
Then 
$F$ satisfies the conditions of Theorem~\ref{tcharsuff}(I) 
with 
$\alpha=\alpha_1$.   
Therefore 
$F\in {\cal A}$.    
\end{thm} 
\proof\
It is easy to see that 
$\supp(F)$ 
must contain 
two consecutive integers 
$k, k+1$, 
and therefore 
$F$ is aperiodic.    
Set 
$\widetilde{c}: = \sum_{j=1}^N c_j \zeta(1+\alpha_j)  >0$, 
$\lambda_j: = \widetilde{c}^{-1} c_j \zeta(1+\alpha_j) 
\in (0,1)$,    
and define   
\[
G: = \sum_{j=1}^N \lambda_j Z_{\alpha_j} \in \PPi(\Zi^+), 
\]
a convex combination of the probabilities 
$Z_{\alpha_j}$. 
It is straightforward to see that 
$G$ satisfies the conditions of 
Theorem~\ref{tcharsuff}(I) with 
$\alpha=\alpha_1$. 
Because  
$F = \widetilde{c} G + P + [F(0)-P(0)]\delta_0$, 
the theorem follows by applying Theorem~\ref{tcharsuff}(II). 
\hfill $\Box$ 

\bigskip

We end this section with a formula connecting 
$A_{\alpha}$ 
with the convolution semigroup 
$(\mu_t^{[\alpha]})_{t\geq 0}$ 
of L\'{e}vy stable probability measures on 
$\Ri^+$. 
For fixed 
$\alpha\in (0,1)$,  
the measures 
$\mu_t^{[\alpha]}$, 
$t\geq 0$,  
are implicitly defined by their Laplace transforms:  
\[
\int_0^{\infty} d\mu_t^{[\alpha]}(s) \, e^{-sz}  = e^{-tz^{\alpha}}
\]
for all 
$z\in \Ci$ with 
$\RRe z \geq 0$ 
(for details see 
\cite[Section IX.11]{Yos}, \cite[Chapter XIII]{Fell2}
and \cite{CaKa}).          
Let 
$(P_s)_{s\geq 0} 
\subseteq \PPi(\Zi^+)$ 
be the Poisson convolution semigroup, 
defined by   
\[
P_s: = e^{-s(\delta_0-\delta_1)} 
= e^{-s} \sum_{k\geq 0} (k!)^{-1} s^k \delta_k.  
\]

\begin{prop} \label{pflevy} 
Let $\alpha\in (0,1)$.   
One has equality 
\begin{equation} 
e^{-t(\delta_0 - A_{\alpha})} 
= \int_0^{\infty} d\mu_t^{[\alpha]}(s)\, P_s    
\label{elevy} 
\end{equation} 
as elements of $\PPi(\Zi^+)$, 
for all 
$t\geq 0$.   
\end{prop} 

Roughly speaking, 
$(e^{-t(\delta_0-A_{\alpha})})_{t\geq 0}$ 
is the convolution semigroup subordinated to the Poisson semigroup 
via  
$(\mu_t^{[\alpha]})_{t\geq 0}$.  
 
\bigskip 

\noindent {\bf Proof of Proposition~\ref{pflevy}.} 
Fix $t>0$.     
Each side of (\ref{elevy}) defines an element of  
$\PPi(\Zi^+)$, and 
it suffices to show that both sides have  
the same generating function.   
Because 
$P_s$ has generating function 
$w\mapsto e^{-s(1-w)}$, 
the right side of (\ref{elevy}) has generating function 
\[
w\mapsto 
\int_0^{\infty} d\mu_t^{[\alpha]}(s) \, e^{-s(1-w)} = e^{-t(1-w)^{\alpha}}. 
\]
From (\ref{eagen}) one sees that 
the left side 
of (\ref{elevy}) also has generating function   
$e^{-t(1-w)^{\alpha}}$.  
\hfill $\Box$

\section{Further examples} \label{s6}

We first construct some examples of elements of 
${\cal A}$ 
which do not fall under the hypotheses of Theorem~\ref{tcharsuff}.  
In these examples 
$1-\phi_F(w)$ 
has a type of logarithmic decay as $w\to 1$ 
within the disc 
$\overline{\Di}$.

\begin{thm} \label{tsubex} 
Fix $\varepsilon\in (0,1]$ and set 
\begin{equation} 
B: = \varepsilon^{-1} \int_0^{\varepsilon} d\alpha \, A_{\alpha} 
\in \PPi(\Zi^+) 
\label{ebdef}
\end{equation} 
(with $A_{\alpha}$ defined by (\ref{efdens})).  
For each $\beta\in (0,1)$ put 
$B_{\beta} : = \sum_{k\geq 0} A_{\beta}(k) B^{(k)}$.   
Then $B_{\beta}\in {\cal A}$, 
and 
\begin{equation} 
\lim_{\xi\to 0} \frac{|\xi|^{\alpha}}{|1-\widehat{B_{\beta}}(\xi)|} 
=0 
\label{edomas} 
\end{equation} 
for each $\alpha\in (0,1)$.  
\end{thm} 

Integration of the estimate (\ref{efderiv}) yields a bound  
$|1-\widehat{F}(\xi)| \leq c' |\xi|^{\alpha}$ 
for all $\xi\in [-\pi, \pi]$;  
thus (\ref{edomas}) shows that the probability 
$F=B_{\beta}$ 
cannot satisfy (\ref{efderiv}) for 
$\alpha\in (0,1)$. 

We point out that (\ref{ebdef}) has a formal similarity to a construction of \cite{Ly} 
(see also \cite[Section 2]{KMOT}),   
where a Ritt operator is produced by integrating certain fractional Volterra operators.

\bigskip 

\noindent {\bf Proof of Theorem~\ref{tsubex}.}  
Clearly (\ref{ebdef}) defines an element 
$B\in \PPi(\Zi^+)$.  
Since $A_{\beta}\in {\cal A}$, 
Proposition \ref{paprop} implies that 
$B_{\beta}\in {\cal A}$.  

Using (\ref{eagen}) one computes  
\begin{equation}  
1 - \phi_B(w) = \varepsilon^{-1} \int_0^{\varepsilon} 
   d\alpha\, (1-w)^{\alpha}    
  = \varepsilon^{-1} [{\mathop{\rm Log}}(1-w)]^{-1} ((1-w)^{\varepsilon}-1)
\label{ebgen} 
\end{equation} 
for $w\in \overline{\Di}$; 
here 
$\mathop{\rm Log} \colon \Ci\backslash (-\infty,0] \to \Ci$ 
is the principal branch of the logarithm, 
with the additional convention  
$({\mathop{\rm Log}}(0))^{-1}=0$.   
Replacing $w$ by $e^{-i\xi}$ in (\ref{ebgen}),  because of the logarithmic 
factor one sees that 
\[
\lim_{\xi\to 0}  \frac{|\xi|^{\alpha}}{|1-\widehat{B}(\xi)|}=0
\]
for each 
$\alpha\in (0,1)$.  
But the definition of $B_{\beta}$  
implies that 
\[
\phi_{B_{\beta}}(w) 
= \sum_{k\geq 0} A_{\beta}(k) (\phi_B(w))^k  
= \phi_{A_{\beta}}(\phi_B(w))
= 1 - (1-\phi_B(w))^{\beta}, 
\]
hence 
$1 - \widehat{B_{\beta}}(\xi) = (1-\widehat{B}(\xi))^{\beta}$, 
and (\ref{edomas}) follows. 
\hfill $\Box$ 

\bigskip 

We conjecture that $B$ defined by (\ref{ebdef}) 
is an element of ${\cal A}$,  for any $\varepsilon\in (0,1]$.  
Evidence for this conjecture is that 
the generating function $\phi_B$ satisfies the sectorial 
condition in (\ref{esectf}), 
as is easily verified from (\ref{ebgen}).    
Actually, $\phi_B$ 
(and similarly 
$\phi_{B_{\beta}}$) 
satisfies not only  (\ref{esectf}) 
but also 
\[
\lim_{w\in \overline{\Di}, w\to 1} {\mathop{\rm Arg}}(1-\phi_B(w))=0. 
\]
This contrasts with the probabilities 
$A_{\alpha}$, $\alpha\in (0,1)$, 
which satisfy 
(use (\ref{eagen})) 
\[
\limsup_{w\in \overline{\Di}, w\to 1} |\mathop{\rm Arg} (1-\phi_{A_{\alpha}}(w))| 
= \alpha\pi/2. 
\]

The construction of Theorem~\ref{tsubex} can be iterated 
to produce further examples of elements of ${\cal A}$.    
For example, 
if for fixed 
$\varepsilon\in (0,1]$, $\beta\in (0,1)$ 
we set  
\[
C: = \varepsilon^{-1} \int_0^{\varepsilon} d\alpha\, B_{\alpha}, 
\;\;\;  
C_{\beta} := \sum_{k\geq 0} A_{\beta}(k) C^{(k)}, 
\] 
then 
$C_{\beta}\in {\cal A}$,   
and 
$1-\phi_{C_{\beta}}(w)= (1-\phi_C(w))^{\beta}$ 
has a type of iterated logarithmic decay 
(slower than logarithmic decay) 
as $w\to 1$. 
We leave details to the reader.  

We next give an example in a quite different direction.

\begin{thm} \label{tinfarg}  
There exists an aperiodic probability $F\in \PPi(\Zi^+)$  
satisfying $\sum_{k\geq 0} kF(k)=+\infty$ 
(that is, 
(\ref{einfmean})) and such that  
\begin{equation} 
\lim_{\xi\to 0} |\mathop{\rm Arg}(1-\widehat{F}(\xi))| = \pi/2. 
\label{earglim} 
\end{equation}  
Thus $F$ does not satisfy (\ref{esectf}), so $F\notin {\cal A}$. 
\end{thm} 
\proof\
Define $F\in L^1(\Zi^+)$ 
by 
$F(0)=F(1)=0$ 
and 
\[
F(k) = \frac{1}{k(k-1)} 
\]
for $k\in \{2, 3, 4, \ldots \}$. 
Clearly 
$\sum_{k\geq 0} k F(k)= +\infty$.   
Using the expansion 
$\mathop{\rm Log}(1-w) = - \sum_{k\geq 1} k^{-1} w^k$ 
one checks that 
\[ 
1- \phi_F(w) 
= 1 - \sum_{k\geq 2} \frac{w^k}{k(k-1)} 
= (1-w)[ 1 - \mathop{\rm Log}(1-w)] 
\]
for all 
$w\in \overline{\Di}\backslash\{1\}$. 
Taking the limit as $w\to 1$ yields 
$\phi_F(1)=1$, 
so 
$F\in \PPi(\Zi^+)$.  
By setting $w=e^{-i\xi}$  in the formula for 
$1-\phi_F(w)$  
one easily verifies (\ref{earglim}).  
\hfill $\Box$ 

\bigskip 

Similarly, it is possible to show 
that the probability 
$Z_1(k): = \zeta(2)^{-1} k^{-2}$, 
$k\in\Ni$ 
(which is the case $\alpha=1$ of (\ref{ezdef})) 
satisfies the statements of Theorem~\ref{tinfarg}.

 \section{Kreiss operators and fractional powers} \label{s7}

 The aim of this section is to establish a generalization of Theorem~\ref{taalpha} 
for Kreiss operators, 
 and then to prove Theorem~\ref{tkritt}.   
Our approach is based on 
the theory of fractional powers of linear operators 
(cf. \cite{MCSA}) 
 and is essentially independent of the theory of previous sections. 
In particular, we do not use the idea of subordination via probabilities.

We need some preliminaries.   
A densely defined, closed  linear operator $V$ in the complex Banach space 
$X$ is said to be of 
{\it type} $\omega$, 
where 
$\omega\in [0, \pi)$,   
if $\sigma(V)\subseteq \overline{\Lambda}_{\omega}$ 
and 
\[
\sup_{\lambda\in {\Lambda}_{\pi-\omega-\varepsilon}} 
\|  \lambda (\lambda I + V)^{-1} \| < \infty 
\]
for every $\varepsilon\in (0, \pi-\omega)$.  
It is a standard fact 
(\cite[Section~2.5]{Dav80}) 
that 
$V$ is of type $\omega$ for some $\omega\in [0, \pi/2)$ 
if and only if 
$-V$ is the generator of a bounded analytic semigroup 
$(e^{-tV})_{t\geq 0}$.

For an operator $V$ of type $\omega$ the fractional powers  
$V^{\alpha}$ can be defined for all $\alpha>0$; 
see \cite{Haa,MCSA} for recent 
expositions of this theory.  
One has relations 
$V^{\alpha} V^{\beta} = V^{\alpha+\beta}$ for $\alpha,\beta>0$, 
and 
$(V^{\alpha})^{\beta} = V^{\alpha \beta}$ for $\alpha\in (0,1)$, $\beta>0$.   
The spectral mapping theorem for fractional powers states that  
\begin{equation} 
\sigma(V^{\alpha}) = \{ z^{\alpha}\colon z\in \sigma(V) \},  
\;\;\; 
\alpha>0.   
\label{efracsp} 
\end{equation} 
If, in addition, the operator $V$ is bounded, 
then $V^{\alpha}$ is bounded for $\alpha>0$.

The next  result
appears for example in \cite[Corollary 3.10]{BBL}, \cite[Proposition 3.1.2]{Haa},  
and in case $\alpha\in (0,1)$ in 
\cite{K5}. 
Since this result is important for our purposes, we sketch a proof.

\begin{thm} \label{ttypefrac} 
Let $V$ be of type $\omega\in (0,\pi)$. 
Then for each $\alpha\in (0, \pi/\omega)$ 
the operator 
$V^{\alpha}$ is of type $\omega \alpha$.  
\end{thm} 
\proof\
By the formula 
$V^{\beta} = (V^{\beta/2^n})^{2^n}$ 
for 
$\beta>0$, $n\in\Ni$ with $\beta<2^n$,  
the proof of the theorem easily reduces to the 
two cases $\alpha\in (0,1)$ 
and $\alpha=2$. 

In case $\alpha\in (0,1)$ 
we use the resolvent formula of Kato \cite{K5} 
\[
(\lambda I + V^{\alpha})^{-1} 
= \frac{\sin \alpha\pi}{\pi} \int_0^{\infty} 
dt \, 
t^{\alpha}
\left[ (\lambda    + e^{i\alpha\pi} t^{\alpha})(\lambda + e^{-i\alpha\pi}t^{\alpha}) 
\right]^{-1}  
(tI + V)^{-1} 
\]
where the integral converges for each 
$\lambda\in \Lambda_{(1-\alpha)\pi}$.   
For  $V$ of  type $\omega$ 
and $|\theta|<\pi-\omega$,   
replace $V$ by $e^{-i\theta}V$ in this formula;  
one can then estimate the norm of the operator 
$(\lambda I + e^{-i \alpha\theta} V^{\alpha})^{-1} 
= e^{i\alpha\theta} (\lambda e^{i\alpha\theta} I + V^{\alpha})^{-1}$ 
for $\lambda\in  \Lambda_{(1-\alpha)\pi}$, 
to see that  
$V^{\alpha}$ is type $\alpha\omega$.  

For 
$V$ of type 
$\omega\in (0,\pi/2)$, 
one can see that $V^2$ is type $2\omega$ 
via the identity 
\[
(te^{i\theta}I  +  V^2)^{-1} 
= (it^{1/2}e^{i\theta/2}I + V)^{-1} (-it^{1/2}e^{i\theta/2}I + V)^{-1}
\]
valid for $t>0$ and $|\theta| < \pi- 2\omega$.  
\hfill $\Box$

\bigskip

\begin{thm} \label{tkreiss} 
Suppose that 
$T\in {\cal L}(X)$ 
is a Kreiss operator.  
Then 
$(I-T)^{\alpha}$ 
is of type $\alpha \pi/2$ 
for all $\alpha \in (0,1]$, 
and  
$I-(I-T)^{\alpha}$ 
is a Ritt operator 
for all $\alpha\in (0,1)$. 
\end{thm}  
\proof\
One has 
$|\mu+1|-1 \geq \RRe \mu$ 
for all $\mu\in \Ci$.  
Therefore the Kreiss resolvent condition implies that   
\[
\| (\mu I + (I-T))^{-1} \| 
\leq c (\RRe \mu)^{-1} 
\]
whenever $\RRe \mu>0$. 
It follows that $I-T$ is type $\pi/2$,  
and hence 
$(I-T)^{\alpha}$ is of type 
$\alpha \pi/2$ for $\alpha\in (0,1)$.  
Thus for $\alpha\in (0,1)$, 
the operator  
$S: = I - (I-T)^{\alpha} \in {\cal L}(X)$   
is such that 
$(e^{-t(I-S)})_{t\geq 0}$ 
is a bounded analytic semigroup.  
Also, from 
$\sigma(T)\subseteq\overline{\Di}$  
and  the spectral mapping theorem (\ref{efracsp})  
one obtains   
$\sigma(S)\subseteq \{1-(1-z)^{\alpha}\colon z\in \overline{\Di}\} 
\subseteq \Di\cup\{1\}$.   
By Theorem~\ref{tritt}, 
$S$ is a Ritt operator. 
\hfill $\Box$  

\bigskip 

\noindent {\bf Remark.}  
The hypothesis in Theorem~\ref{tkreiss} that $T$ be Kreiss can be considerably 
weakened: 
the conclusions of the theorem hold whenever 
$T\in {\cal L}(X)$ 
is such that $I-T$ is of type $\pi/2$ and 
$\sigma(T)\subseteq \overline{\Di}$.  

\bigskip

The next result is a corollary of Theorem~\ref{ttypefrac}.

\begin{cor} \label{csgpp} 
If $-V$ is the generator of a bounded analytic semigroup 
$(e^{-tV})_{t\geq 0}$, 
then there exists $\gamma_0>1$ 
such that $-V^{\gamma}$ 
is the generator of a bounded analytic semigroup 
$(e^{-tV^{\gamma}})_{t\geq 0}$ 
for each 
$\gamma\in (1, \gamma_0)$. 
\end{cor}    
\proof\ 
Because $V$ is of type $\omega$ for some 
$\omega \in (0, \pi/2)$, 
then 
$V^{\gamma}$ is of type $\omega \gamma < \pi/2$ 
whenever  
$\gamma \in (1, \pi/(2\omega))$. 
\hfill $\Box$ 

\bigskip

\noindent {\bf Proof of Theorem \ref{tkritt}.}    
(IV)$\Rightarrow$(III) is proved by 
setting $S: = I - (I-T)^{\gamma}$ 
for some $\gamma\in (1, \gamma_0)$, 
noting that 
$I-T = (I-S)^{1/\gamma}$.  
(III)$\Rightarrow$(II) 
is immediate since every power-bounded operator 
is a Kreiss operator.  
(II)$\Rightarrow$(I) is just Theorem~\ref{tkreiss}. 

Let us prove (I)$\Rightarrow$(IV). 
Because $T$ is Ritt 
the semigroup $(e^{-t(I-T)})_{t\geq 0}$ 
is bounded analytic, 
and Corollary \ref{csgpp} implies that 
the semigroup $(e^{-t(I-T)^{\gamma}})_{t\geq 0}$ 
is bounded analytic 
for $\gamma>1$ sufficiently close to $1$.    
Also, by (\ref{efracsp}) and by  
the last statement of Theorem~\ref{tritt} applied to 
$\sigma(T)$, 
it is easy to see that 
\[
\sigma(I-(I-T)^{\gamma}) = \{1-(1-z)^{\gamma}\colon z\in \sigma(T) \} 
\subseteq \Di\cup\{1\}
\] 
for $\gamma>1$ sufficiently close to $1$.  
Thus by Theorem~\ref{tritt}, the operator 
$I-(I-T)^{\gamma}$ is a Ritt operator for 
$\gamma>1$ 
sufficiently close to $1$. 
\hfill $\Box$

\bigskip 

Finally, Corollary~\ref{csingle} is a consequence of 
Theorem~\ref{tkreiss}, the spectral mapping 
(\ref{efracsp}), 
and the observation that   
$(I-T)^{\alpha}\neq 0$     
whenever $\alpha\in (0,1)$ and 
$T\neq I$.   
(For, if $(I-T)^{\alpha}=0$ 
then 
$I-T = ((I-T)^{\alpha})^{1/\alpha} = 0$.)

\section{Appendix} \label{sapp}    

In this appendix we give a proof of Proposition~\ref{papchar}. 
Throughout we consider an adapted probability   
$F\in \PPi(\Zi)$.  
Observe the properties    
$\widehat{F}(0)=1$ 
and 
$|\widehat{F}(\xi)|\leq 1$ 
for all 
$\xi\in [-\pi, \pi]$.

\bigskip

(I)$\Rightarrow$(II):   
suppose that $F$ is aperiodic and
$\eta\in [-\pi, \pi]$ 
with  
$|\widehat{F}(\eta)|=1$; 
we must show that $\eta=0$. 
Setting 
$\tau: =\widehat{F}(\eta) = 
\sum_{k\in \Zi} F(k) e^{-ik\eta}$, 
since 
$|\tau|=1$ it follows that  
\[
e^{-ik\eta} = \tau
\]
for all 
$k\in \supp(F)$.   
After fixing a $k'\in \supp(F)$, 
then  
$e^{-i(k-k')\eta}=1$ 
so that 
$\eta (k-k') \in 2\pi \Zi := \{ 2\pi m\colon m\in \Zi\}$ 
for all 
$k \in \supp(F)$. 
Aperiodicity implies that the set 
$\{k-k'\colon k\in \supp(F) \}$ 
generates the group 
$\Zi$. 
Therefore $\eta \Zi \subseteq 2\pi \Zi$, 
so that 
$\eta\in 2\pi\Zi$   
and  
$\eta=0$.

(II)$\Rightarrow$(III) 
is trivial, 
since 
$\widehat{F}(0)=1$.

(III)$\Rightarrow$(I): 
if $F$ is not aperiodic then 
there exist 
$m\in \{2, 3, 4, \ldots \}$   
and 
$r\in \{0, 1, \ldots, m-1\}$ 
such that 
$\supp(F) \subseteq m\Zi + r$. 
But 
$r\neq 0$   
since 
$F$ is adapted.    
Then  
\[
\widehat{F}(2\pi/m) = \sum_{l\in \Zi} F(ml+r) e^{-i(ml+r)2\pi/m} 
= e^{-2\pi ir/m}  \notin \Di \cup \{1\},  
\]
and condition (III) fails.     

(IV)$\Rightarrow$(III): 
power-boundedness of $G$ implies 
(recall (\ref{especteq})) 
that 
$\{ \widehat{G}(\xi) \colon \xi\in [-\pi, \pi] \}
\subseteq \overline{\Di}$. 
Then    
$\widehat{F}(\xi) = \beta \widehat{G}(\xi) + (1-\beta) 
\in \Di\cup \{1\}$ 
for all 
$\xi\in [-\pi, \pi]$.  
   
\bigskip

To prove 
(I)$\Rightarrow$(IV),  
we begin with the special case where 
$F$ has a finite second moment.

\begin{lemma} \label{lsecmom} 
If 
$F\in \PPi(\Zi)$ 
is aperiodic and 
$\sum_{k\in \Zi} k^2 F(k)< \infty$, 
then $F$ satisfies condition (IV) 
of Proposition~\ref{papchar}.  
\end{lemma} 
\noindent {\bf Proof of Lemma~\ref{lsecmom}.}   
Set 
$a: = \sum_{k\in\Zi} k F(k) \in \Ri$. 
For each 
$\beta>0$  
define $F_{\beta}\in L^1(\Zi)$ 
and $\tau_{\beta}\in C([-\pi, \pi])$ 
by 
\[
F_{\beta}: = \beta^{-1} ( F- (1-\beta)\delta_0), 
\;\;\; 
\tau_{\beta}(\xi): =  e^{i\beta^{-1} a\xi} \widehat{F}_{\beta}(\xi)
 = \sum_{k\in \Zi} F_{\beta}(k) e^{-i(k-\beta^{-1}a)\xi}  
\]
for $\xi\in [-\pi, \pi]$. 
Note that $F_1=F$.     
To obtain the lemma we must show that $F_{\beta}$ is power-bounded 
for some $\beta\in (0,1)$.  

We first show
that there exist  
$\beta_0 \in (0,1)$ and 
$b>0$ such that  
\begin{equation} 
|\tau_{\beta}(\xi)| = |\widehat{F}_{\beta}(\xi)| \leq 1 - b\xi^2 
\label{ebetaest} 
\end{equation} 
for all 
$\xi\in [-\pi, \pi]$ 
and 
$\beta\in (\beta_0, \beta_0^{-1})$.  
Because  
$\sum_{k\in \Zi}  
k^2 F(k)< \infty$, 
one easily checks that   
$\tau_{\beta} \in C^2([-\pi, \pi])$ 
and that   
$\tau_{\beta}(0)=1$,   
$\tau_{\beta}'(0)=0$ 
for all 
$\beta>0$.   
Thus Taylor's formula yields 
\begin{equation} 
\tau_{\beta}(\xi) = 1 + \xi^2 \int_0^1 ds\, \tau_{\beta}''(s\xi)(1-s) 
\label{etayl} 
\end{equation} 
for all 
$\xi\in [-\pi, \pi]$ and 
$\beta>0$.   
Observe that the functions 
$\tau_{\beta}''$ converge uniformly to 
$\tau_1''$ 
as 
$\beta\to 1$, 
and that 
$\tau_1''(0) = - \sum_{k\in \Zi} (k-a)^2 F(k) < 0$.   
Then one deduces from (\ref{etayl}) 
that there is a $\delta\in (0,1)$ 
such that an estimate (\ref{ebetaest}) 
holds for all $\xi\in [-\delta, \delta]$ 
when  $\beta$ is sufficiently close to $1$. 
Aperiodicity of $F$ implies that 
\[
\sup \{ |\widehat{F}(\xi)| \colon \xi\in [-\pi, -\delta] \cup [\delta, \pi] \} 
< 1,     
\] 
so it is also easy to obtain an estimate of type (\ref{ebetaest}) 
when $\xi\in [-\pi, -\delta]\cup [\delta,\pi]$ 
and $\beta$ is close to $1$. 
Thus (\ref{ebetaest}) is obtained for some $\beta_0\in (0,1)$.

For the rest of the proof, 
fix 
$\beta\in (\beta_0, 1)$ and put
$G:= F_{\beta}$.  
The Plancherel formula gives   
$\sum_{k\in \Zi}  |H(k)|^2 
= (2\pi)^{-1} \int_{-\pi}^{\pi} d\xi\, |\widehat{H}(\xi)|^2$ 
for all 
$H\in L^1(\Zi)$.
Using (\ref{ebetaest}) and the inequality 
$1-b\xi^2 \leq e^{-b \xi^2}$, 
we obtain an estimate  
\begin{eqnarray} 
\sum_{k\in \Zi}  |G^{(n)}(k)|^2 
& = & (2\pi)^{-1} \int_{-\pi}^{\pi} d\xi\, 
|(\widehat{F}_{\beta}(\xi))^n|^2  \nonumber  \\
& \leq & (2\pi)^{-1} \int_{-\pi}^{\pi} d\xi\, e^{-2bn\xi^2}  
\leq c n^{-1/2} 
\label{egest1} 
\end{eqnarray} 
for all 
$n\in\Ni$.  
Writing  
$\partial_{\xi}$ 
for differentiation with respect to 
$\xi$, 
observe that the 
$\Zi$-Fourier transform of 
$k\mapsto -i(k-\beta^{-1} an) G^{(n)}(k)$ 
is 
\begin{eqnarray*} 
(\partial_{\xi} + i\beta^{-1} an) ( \widehat{F}_{\beta}(\xi)^n) 
 & = & e^{-i\beta^{-1} an\xi} \partial_{\xi} ( \tau_{\beta}(\xi)^n)     \\
 & = & e^{-i\beta^{-1} an\xi} n  \tau_{\beta}(\xi)^{n-1} \tau_{\beta}'(\xi).     
\end{eqnarray*} 
Because 
$\tau_{\beta}$ is $C^2$ and 
$\tau_{\beta}'(0)=0$,  
there is an estimate  
\[
|\tau_{\beta}'(\xi)| = |\tau_{\beta}'(\xi) - \tau_{\beta}'(0)| 
\leq c' |\xi|, 
\;\;\; 
\xi\in [-\pi, \pi].   
\]     
By these facts and (\ref{ebetaest}),   
\begin{eqnarray} 
\sum_{k\in \Zi} (k-\beta^{-1} an)^2 |G^{(n)}(k)|^2 
 & = & (2\pi)^{-1} \int_{-\pi}^{\pi} d\xi\, 
| n \tau_{\beta}(\xi)^{n-1} \tau_{\beta}'(\xi) |^2    \nonumber \\
& \leq  & c \int_{-\pi}^{\pi} d\xi\, 
 n^2 \xi^2 e^{-2b(n-1)\xi^2}    
\leq c'  n^{1/2} 
\label{egest2} 
\end{eqnarray} 
for all 
$n\in\Ni$. 
The Cauchy-Schwarz inequality, (\ref{egest1}) and (\ref{egest2}) 
imply that 
\begin{eqnarray*} 
& & \sum_{k\in \Zi} |G^{(n)}(k)|   \\ 
& \leq & 
\left( \sum_k [1+n^{-1} (k-\beta^{-1} an)^2]^{-1} \right)^{1/2} 
\left( \sum_k [1+n^{-1} (k-\beta^{-1} an)^2] |G^{(n)}(k)|^2 \right)^{1/2}  \\
& \leq & c n^{1/4} 
 \left( \sum_k [1+n^{-1} (k-\beta^{-1} an)^2] |G^{(n)}(k)|^2 \right)^{1/2}  
  \leq c' 
\end{eqnarray*} 
for all $n\in\Ni$. 
Thus $G$ is power-bounded.    
\hfill $\Box$ 
  
\bigskip 

We next prove 
condition (IV) 
in case $F$ is aperiodic 
and 
$\sum_{k\in \Zi} k^2 F(k) = \infty$. 
In this case set 
$\gamma: = \sum_{k\in \Zi}  (1+k^2)^{-1} F(k) \in (0,1)$,  
and define 
$\widetilde{F}\in \PPi(\Zi)$ 
by 
\[
\widetilde{F}(k) = \gamma^{-1} (1+k^2)^{-1} F(k), 
\;\;\; 
k\in \Zi. 
\]
Since 
$\widetilde{F}$ is aperiodic  
and 
$\sum_{k\in\Zi} k^2 \widetilde{F}(k) < \infty$, 
by Lemma~\ref{lsecmom} there exist  
$\beta\in (0,1)$ 
and a power-bounded 
$\widetilde{G}\in L^1(\Zi)$ 
such that  
\[
\widetilde{F} = \beta \widetilde{G} + (1-\beta) \delta_0.   
\]
Because  
$F(k)\geq \gamma \widetilde{F}(k)$ for all 
$k\in \Zi$, 
we have  
$F = \gamma \widetilde{F} + (1-\gamma) H$ 
for a certain 
$H \in \PPi(\Zi)$. 
Setting   
$\alpha: = 1-\gamma +\beta\gamma \in (0,1)$ 
and 
\[
G: = \frac{\beta\gamma}{\alpha} \widetilde{G} + \frac{1-\gamma}{\alpha} H 
\in L^1(\Zi),  
\]
it follows that   
$F = \alpha G + (1-\alpha) \delta_0$. 
Since $G$ is a convex combination of the power-bounded elements 
$\widetilde{G}$, 
$H$ 
in the (commutative) Banach algebra $L^1(\Zi)$, 
it is easy to see that 
$G$ is power-bounded.   
Thus 
$F$ satisfies condition (IV).   

To complete the proof of Proposition~\ref{papchar}, 
we verify (\ref{eaphalf}).    
Let $F$ be aperiodic,  
let $G$ be as in condition (IV),     
and apply Theorem~\ref{tnevhalf} with the operators 
$S: = L(F)$ and $T: = L(G)$ acting in $X=L^1(\Zi)$.    
Then (\ref{eaphalf}) follows. 

\bigskip 

\noindent {\bf Remark.}     
One can prove (\ref{eaphalf}) directly   
without using Theorem~\ref{tnevhalf}, 
at least in the case  
$\sum_{k\in\Zi} k^2 F(k)< \infty$.  
For in that case, 
by arguing in a similar way to the proof of Lemma~\ref{lsecmom} 
one can show that
\[
\sum_{k\in\Zi} [1 + n^{-1}(k-an)^2]  | F^{(n)}(k) - F^{(n+1)}(k) |^2  
\leq c n^{-3/2} 
\]
for all 
$n\in\Ni$, 
where 
$a: = \sum_{k\in\Zi} kF(k)$.  
Then (\ref{eaphalf}) follows by the Cauchy-Schwarz inequality 
as in the proof of Lemma~\ref{lsecmom}.

\bigskip 

Department of Mathematics 

Macquarie University 

NSW 2109

Australia

Email: ndungey@ics.mq.edu.au

\end{document}